\documentclass[12pt]{amsart}
\usepackage{amssymb}
\usepackage{graphics}

\numberwithin{equation}{section}

\def\bea{\begin{align}}
\def\eea{\end{align}}

\newtheorem{Thm}{Theorem}[section]
\newtheorem{Prop}[Thm]{Proposition}
\newtheorem{Lem}[Thm]{Lemma}
\newtheorem{Cor}[Thm]{Corollary}

\theoremstyle{remark}
\newtheorem{Rem}[Thm]{Remark}
\newtheorem*{Ack}{Acknowledgment}

\theoremstyle{definition}
\newtheorem{Def}[Thm]{Definition}

\newtheorem{Conv}[Thm]{Convention}
\newtheorem{Conj}[Thm]{Conjecture}

\def\bconv{\begin{Conv}}
\def\econv{\end{Conv}}
\def\bprop{\begin{Prop}}
\def\eprop{\end{Prop}}
\def\bconj{\begin{Conj}}
\def\econj{\end{Conj}}
\def\bth{\begin{Thm}}
\def\eth{\end{Thm}}
\def\bcor{\begin{Cor}}
\def\ecor{\end{Cor}}
\def\blem{\begin{Lem}}
\def\elem{\end{Lem}}
\def\bdefiniz{\begin{Def}}
\def\edefiniz{\end{Def}}

\newcommand{\dd}{\mathrm{d}}

\newcommand{\lora}{{\longrightarrow}}

\newcommand{\R}{\mathbb{R}}
\newcommand{\reali}{\mathbb{R}}

\newcommand{\N}{\mathbb{N}}

\newcommand{\imbr}[1]{{\,\hbox{\rm Imb}\,(S^1,\reali^{#1}) }}
\newcommand{\imbrr}[1]{{\,\hbox{\rm Imb}\,(\reali,\reali^{#1}) }}

\def\cmp#1,{{ Commun.\ Math.\ Phys.\ \bf #1},}
\def\jmp#1,{{ J.\ Math.\ Phys.\ \bf #1},}
\def\pl#1,{{ Phys.\ Lett.\ \bf #1},}
\def\npb#1,{{ Nucl.\ Phys.\ {\bf B #1}},}
\def\mpl#1,{{ Mod.\ Phys.\ Lett.\ \bf #1},}
\def\pr#1,{{ Phys.\ Rev.\ \bf #1},}
\def\prl#1,{{ Phys.\ Rev.\ Lett.\ \bf #1},}
\def\lmp#1,{{ Lett.\ Math.\ Phys.\ \bf #1},}
\def\jktr#1,{{ Jour.\ of Knot Theory and its Ramification\ \bf #1},}
\def\bams#1,{{ Bull.\ Amer.\ Math.\ Soc.\ \bf #1},}
\def\agt#1,{{Algebr.\ Geom.\ Topol. \bf #1},}

\begin{document}

\title[Algebraic structures on graph cohomology]
{Algebraic structures on graph cohomology}

\author[A.~S.~Cattaneo]{Alberto S.~Cattaneo}
\address{Mathematisches Institut --- Universit\"at Z\"urich--Irchel ---
Winterthurerstrasse 190 --- CH-8057 Z\"urich --- Switzerland}
\email{asc@math.unizh.ch}

\author[P.~Cotta-Ramusino]{Paolo Cotta-Ramusino}
\address{Dipartimento di Fisica --- Universit\`a degli Studi di
Milano \& INFN Sezione di Milano ---  Via Celoria, 16 --- I-20133 Milano 
--- Italy}
\email{cotta@mi.infn.it}

\author[R. Longoni]{Riccardo~Longoni}
\address{Dipartimento di Matematica ``G. Castelnuovo'' --- 
Universit\`a di Roma ``La Sapienza'' ---  Piazzale Aldo Moro, 2 
--- I-00185 Roma --- Italy}
\email{longoni@mat.uniroma1.it}

\thanks{A.~S.~C. acknowledges partial support of SNF grant
  No.~20-100029/ 1}


\begin{abstract}
We define algebraic structures on graph cohomology and prove that they
correspond to algebraic structures on the cohomology 
of the spaces of imbeddings of $S^1$ or $\R$ into $\R^n$. As a corollary, 
we deduce the existence of an infinite number of nontrivial cohomology 
classes in $\imbr n$ when $n$ is even and greater than 3. Finally, we 
give a new interpretation of the anomaly term for the Vassiliev invariants
in $\R^3$. 
\end{abstract}


\maketitle

\section{Introduction}

In this paper we consider the spaces $\imbr n$ of imbeddings of
the circle into $\R^n$ and the spaces $\imbrr n$ of imbeddings of
the real line into $\R^n$ with fixed behavior at infinity, namely,
imbeddings that coincide with a fixed imbedded line in $\R^n$
outside a compact subset.  

When $n>3$ one can define, using configuration space integrals, 
chain maps from certain graph complexes $(\mathcal D_n^{k,m},\delta)$ 
or $(\mathcal L_n^{k,m},\delta)$ to the de~Rham 
complexes of the spaces of imbeddings:
\begin{align}
\label{chainmap1}
I\colon(\mathcal D_n^{k,m},\delta) &\to (\Omega^{(n-3)k+m} 
(\imbr n), \dd),\\
\label{chainmap2}
I\colon(\mathcal L_n^{k,m},\delta) &\to (\Omega^{(n-3)k+m} 
(\imbrr n), \dd).
\end{align}
We will use the same symbol $I$ also for the induced maps in cohomology. 
Our main interest is to determine which algebraic structures are 
preserved by the maps $I$. 

On the one side, in fact, the de~Rham complexes are differential graded 
commutative algebras. Moreover, there exist a multiplication $\imbrr n \times 
\imbrr n \to \imbrr n$ given by attaching the end of the first imbedding to 
the beginning of the second, and rescaling.
This operation is non-associative but gives rise, together with the 
cup product, to a Hopf algebra structure on $H(\imbrr n)$ for $n>3$. 
Indeed, as shown in \cite{B}, the multiplication is part of an action of 
the operad of little 2-cubes on $\imbrr n$.

On the other side, we show in Section~\ref{sec:alg} that
$(\mathcal D_n^{k,m},\delta)$ and $(\mathcal L_n^{k,m},\delta)$ 
(whose definitions are recalled in Section~\ref{grcoho}) can be endowed, 
respectively, with a differential algebra and a differential 
Hopf algebra structure. The existence of these algebraic structures 
implies in particular the existence of infinitely many nontrivial 
classes of trivalent cocycles of even type (see Section~\ref{sec:inf}), 
and, as a consequence of the results of \cite{CCRL}, of infinitely 
many nontrivial elements of 
$H(\imbr n)$, for every even $n\ge 4$.

The central result of the paper, contained in 
Section~\ref{sec:hopf}, is that, for $n>3$, the maps $I$ in cohomology 
respect all the above algebraic structures. 

Finally, in Section~\ref{sec:anomaly} we consider the case $n=3$ and 
$m=0$ (the case of Vassiliev invariants) and we give a new interpretation 
of the so-called ``anomaly'' term \cite{BT, AF} as the obstruction for the
Bott--Taubes map from trivalent cocycles to Vassiliev invariants 
to be a coalgebra map.

As a final remark, we recall that there exist other approaches to study 
the cohomology of spaces of imbeddings based on graph cohomology 
\cite{Si,T}. It would be interesting to understand our results in these 
other contexts.

\begin{Ack}
We thank Domenico Fiorenza, Dev Sinha, Jim Stasheff, Victor Tourtchine 
and the Referee for useful comments and suggestions on a first version 
of this paper. 
A.~S.~C. acknowledges the University of Roma ``La Sapienza'' (special 
trimester ``Moduli, Lie theory, interactions with physics''), R.~L. 
and P.~C.-R. acknowledge the University of Zurich, R.~L. acknowledges the 
University of Milano for their kind hospitality during the preparation 
of the work. 
\end{Ack}

\section{Graph cohomology}
\label{grcoho}

We briefly recall the definition of the {\em graph complexes} given in 
\cite{CCRL}, Section 4. A graph consists of an {\em oriented  circle}\/ 
and many 
{\em edges}\/ joining vertices. The vertices lying on the circle are 
called {\em external vertices}, those lying off the circle are called 
{\em internal vertices} and are required to be at least trivalent. We define 
the order $k\ge 0$ of a graph to be minus its Euler characteristic, and the 
degree $m\ge 0$ to be the deviation of the graph from being trivalent, namely:
\begin{align*}
k&=e -v_i,\\
m&=2e-v_e-3v_i,
\end{align*}
where $e$ is the number of edges, $v_i$ the number of internal 
vertices and $v_e$ the number of external vertices of the graph.

The graph complexes $(\mathcal D_n^{k,m},\delta)$ depend only 
on the parity of $n$, and we will denote by $\mathcal D_o^{k,m}$ and 
$\mathcal D_e^{k,m}$ the real vector spaces generated by decorated 
graphs of order $k$ and degree $m$, of odd and even type, respectively.  

The type of a graph depends on the decoration which we put on it.
By definition, the decoration in $\mathcal D_o^{k,m}$ is given by 
numbering all the vertices (with the convention that we first number 
the external vertices and then the internal ones) and orienting the edges. 
Then one takes the quotient by the following relations: a
cyclic permutations of the external vertices or a permutation of the 
internal vertices multiplies the graph by the sign of the permutation, 
a reversal of an orientation of an edge produces a minus sign.
An extra decoration is needed on edges connecting the same external 
vertex, namely an ordering of the two half-edges forming them. The 
decoration in $\mathcal D_e^{k,m}$ is given by numbering the 
external vertices and numbering the edges, while the relations are as 
follows: a cyclic permutations of the external vertices or a 
permutation of the edges multiplies the graph by the sign of the 
permutation. 

Double lines and internal loops are not allowed, namely, 
in $\mathcal D_o^{k,m}$ and $\mathcal D_e^{k,m}$ we quotient 
by the subspace generated by all the graphs containing 
two edges joining the same pair of vertices and by all the 
graphs containing edges whose end-points are the same internal vertex.

By an {\em arc}\/ we mean a piece of the oriented circle between
two consecutive vertices, and by a {\em regular edge}\/ we mean an edge 
with at least one internal end-point.

The coboundary operators $\delta_o\colon \mathcal
D^{k,m}_o \to \mathcal D^{k,m+1}_o$ and $\delta_e\colon \mathcal
D^{k,m}_e \to \mathcal D^{k,m+1}_e$ are linear operators, whose action
on a graph $\Gamma$ is given by the signed sum of all the graphs obtained
from $\Gamma$ by contracting, one at a time, all the regular edges and 
arcs of the graphs. The signs are as follows: if we contract
an arc or edge connecting the vertex $i$ with the vertex $j$, and oriented
from $i$ to $j$, then the sign is $(-1)^{j}$ if $j>i$ or $(-1)^{i+1}$
if $j<i$. If we contract an edge labeled by $\alpha$, then the
sign is $(-1)^{\alpha+1+v_e}$ where $v_e$ is the number of external
vertices of the graph. The decoration on $\Gamma$ induces a decoration on
the contracted graphs as follows: contraction of the edge between the vertex 
$i$ and the vertex $j$ produces a new vertex labeled by $\min(i,j)$ while
the vertex labels greater than $\max(i,j)$ are rescaled by one. Similarly 
when the edge labeled by $\alpha$ is collapsed, the edge labels greater 
than $\alpha$ are rescaled by one.

The complexes $(\mathcal L_o^{k,m},\delta_o)$ and $(\mathcal L_e^{k,m}, 
\delta_e)$ are defined in the same manner, except that instead of graphs 
with an oriented circle we consider graphs with an {\em oriented line}. 

The cohomology groups with respect to the above coboundary operators will be 
denoted by $H^{k,m}(\mathcal D_o)$, $H^{k,m}(\mathcal D_e)$, 
$H^{k,m}(\mathcal L_o)$ and $H^{k,m}(\mathcal L_e)$.
When we do not want to specify the parity of these 
complexes we simply write $\mathcal D^{k,m}$ or $\mathcal L^{k,m}$. 
Moreover we set $\mathcal D = \bigoplus_{k,m}\mathcal D^{k,m}$, 
$\mathcal L = \bigoplus_{k,m}\mathcal L^{k,m}$,
$H(\mathcal D) = \bigoplus_{k,m} H^{k,m}(\mathcal D)$ and 
$H(\mathcal L) = \bigoplus_{k,m} H^{k,m}(\mathcal L)$

\section{Algebraic structures on graphs}
\label{sec:alg}

\subsection{Operations on $\mathcal D^{k,m}$}

Suppose that $\Gamma_1$ and $\Gamma_2$ are two graphs with $v_e(\Gamma_1)$ 
and $v_e(\Gamma_2)$ external vertices, respectively. The sets $V_e^1$ and 
$V_e^2$ of external vertices of $\Gamma_1$ and $\Gamma_2$ respectively,
are cyclically ordered. A $(V_e^1, V_e^2)$-shuffle is a permutation of 
the set $V_e^1\cup V_e^2$ respecting the cyclic order of $V_e^1$ and 
$V_e^2$. If $\sigma$ is a $(V_e^1, V_e^2)$-shuffle, then we can combine 
$\Gamma_1$ and $\Gamma_2$ in a single graph $\Gamma_1\bullet_\sigma\Gamma_2$ 
with $v_e(\Gamma_1)+v_e(\Gamma_2)$ external vertices by placing the legs of 
$\Gamma_1$ and $\Gamma_2$ into distinct external vertices on an oriented 
circle, according to $\sigma$.

We put the labels in $\Gamma_1 \bullet_\sigma \Gamma_2$ as follows: 
the vertices and edges which come from $\Gamma_1$ are
numbered in the same way as in $\Gamma_1$. Next we number the 
vertices coming from $\Gamma_2$ by adding to the corresponding label of
$\Gamma_2$ the number of labeled vertices of $\Gamma_1$ and the edges by 
adding to the corresponding label of $\Gamma_2$ the number of edges of
$\Gamma_1$. If a vertex or edge is not labeled in $\Gamma_1$ or
$\Gamma_2$, it remains unlabeled in $\Gamma_1 \bullet_\sigma \Gamma_2$
%
%
The multiplication $\Gamma_1 \bullet \Gamma_2$ of two graphs
$\Gamma_1$ and $\Gamma_2$ is then defined as
$$\Gamma_1 \bullet \Gamma_2 = (-1)^{\lambda(\Gamma_1,\Gamma_2)}  \sum_\sigma \Gamma_1\bullet_\sigma\Gamma_2$$ where 
the sum is taken over all the $(V_e^1, V_e^2)$-shuffles and 
\begin{equation}
\label{segno}
\lambda(\Gamma_1,\Gamma_2) =
\begin{cases}
v_e(\Gamma_2) \, e(\Gamma_1) & \text{for}\,\,\Gamma_i\in\mathcal D_e\\
0  & \text{for}\,\,\Gamma_i\in\mathcal D_o.
\end{cases}
\end{equation}
In the above formula, $e(\Gamma_1)$ denotes the number of edges 
of $\Gamma_1$.

Extending this product to $\mathcal D$ by linearity, we obtain an
associative operation called the {\em shuffle product}:
\[
\bullet \colon \mathcal D^{k_1,m_1} \otimes \mathcal D^{k_2,m_2}
\to \mathcal D^{k_1+k_2,m_1+m_2}.
\]

We define a new grading $|\cdot|$ of the graphs generating
$\mathcal D^{k,m}_e$ and $\mathcal D^{k,m}_o$, by considering the
total number of labels of the graph: 
\[
|\Gamma|\equiv
\begin{cases}
e(\Gamma) + v_e(\Gamma)&
\text{for}\,\,\Gamma_i\in\mathcal D_e\\ 
v(\Gamma)& \text{for}\,\,\Gamma_i\in\mathcal D_o.
\end{cases}
\]
where $v(\Gamma)$ is the total number of vertices of $\Gamma$. Modulo
2, the new grading is equal to $k+m$ in the even case and to $m$ in
the odd case. With respect to the new grading the integration maps $I$
will be grading-preserving modulo 2. The shuffle product is graded
commutative with respect to this new grading. In fact one can easily see that
$\Gamma_1 \bullet \Gamma_2 = (-1)^{|\Gamma_1|\, |\Gamma_2|}\Gamma_1
\bullet \Gamma_2$. Leibnitz rule also holds between the coboundary
operator $\delta$, which has degree $-1$ with respect to the new grading, 
and the shuffle product, which has degree 0. More precisely 

\bprop
\label{prop:db}
$\delta(\Gamma_1 \bullet \Gamma_2)= \delta(\Gamma_1) \bullet \Gamma_2
+(-1)^{|\Gamma_1|} \Gamma_1 \bullet \delta(\Gamma_2)$ for
every $\Gamma_i\in\mathcal D$.
\eprop

\begin{proof}
Let $\sigma$ be a $(V_e^1, V_e^2)$-shuffle. When we apply the coboundary 
operator $\delta$ to one of the graphs of $\Gamma_1 \bullet_\sigma \Gamma_2$
we can either collapse an edge of $\Gamma_1$, or collapse two external 
vertices of $\Gamma_1$, or collapse an edge of $\Gamma_2$, or collapse 
two external vertices of $\Gamma_2$, or collapse an external vertex of 
$\Gamma_1$ with an external vertex of $\Gamma_2$. The first and second 
contribution yield $\delta(\Gamma_1) \bullet_\sigma \Gamma_2$, while the 
third and fourth yield $\Gamma_1 \bullet_\sigma \delta(\Gamma_2)$. Taking 
into account the signs we obtain the formula
\[
\delta(\Gamma_1 \bullet_\sigma \Gamma_2) = \delta(\Gamma_1) \bullet_\sigma 
\Gamma_2 + (-1)^{|\Gamma_1|} \Gamma_1 \bullet_\sigma \delta(\Gamma_2) +
R_\sigma(\Gamma_1,\Gamma_2)
\]
where $R_\sigma(\Gamma_1,\Gamma_2)$ is a linear combination of the graphs 
obtained from $\Gamma_1 \bullet_\sigma \Gamma_2$ by contracting an 
external vertex of $\Gamma_1$ with an external vertex of $\Gamma_2$.
%
%
Now, suppose that $\sigma$ brings the $i$th external vertex of $\Gamma_1$
next to the $j$th external vertex of $\Gamma_2$ in
$\Gamma_1\bullet_\sigma\Gamma_2$. Collapsing the arc between these two 
vertices yields a contribution, called $\Lambda$, to
$R_{\sigma}(\Gamma_1,\Gamma_2)$. Then consider the $(V_e^1, V_e^2)$-shuffle 
$\sigma^\tau$ obtained by composing $\sigma$ with a transposition of the 
two vertices considered above. Collapsing the arc in 
$\Gamma_1\bullet_{\sigma^\tau}\Gamma_2$ between these vertices yields the 
contribution $-\Lambda$ to $R_{\sigma}(\Gamma_1,\Gamma_2)$. This fact in 
turn implies that $\sum_\sigma R_\sigma(\Gamma_1,\Gamma_2)=0$, 
with the sum taken over all $(V_e^1, V_e^2)$-shuffles.
\end{proof}

As a consequence, $(\mathcal D, \delta, \bullet)$ is a differential 
graded commutative algebra. It also has a unit $\mathbf 1$ consisting of
the graph without edges. It follows that its cohomology $H(\mathcal D)$
is a graded commutative algebra with unit.

\subsection{Operations on $\mathcal L^{k,m}$}

The product $\bullet$ on $\mathcal L^{k,m}$ is defined exactly as in
the previous case, i.e., as the shuffle product of two graphs.
The unit {\bf 1} is the graph with no edges. In addition we have a 
comultiplication $\Delta$ mapping a graph $\Gamma$ to the signed
sum of $\Gamma'\otimes \Gamma''$ over all possible ways to cut $\Gamma$ 
into two disconnected parts $\Gamma'$ and $\Gamma''$ by removing an 
internal point of one of the arcs of $\Gamma$. 
By convention, if the oriented line is oriented, say, from left to 
right, then $\Gamma_1$ is the left-most component of 
$\Gamma\setminus\{\mbox{pt}\}$. A graph $\Gamma$ is primitive if it cannot be 
disconnected in a nontrivial way 
by removing an internal point of an arc, and in this case 
$\Delta\Gamma = \mathbf{1}\otimes \Gamma + \Gamma \otimes \mathbf{1}$. 
To fix the signs in the general case, first of all we order the primitive 
subgraphs of $\Gamma$ with respect to the oriented line, and assume that 
the labels are compatible with this ordering (i.e., all labels of a 
primitive subgraph are less than any label of a subsequent one). 
The decoration induced by $\Gamma$ on $\Gamma'$ and $\Gamma''$ is 
determined by rescaling the labels of $\Gamma''$ and 
\[
\Delta\Gamma = \sum (-1)^{\lambda(\Gamma',\Gamma'')} 
\Gamma'\otimes \Gamma'',
\]
where $\lambda(\Gamma',\Gamma'')$ is given in eq.~\eqref{segno}.

We also define a counit $\epsilon$ by $\epsilon(\Gamma)=0$ if 
$\Gamma\neq \mathbf 1$ and $\epsilon(\mathbf 1)=1$.

\bth
With the above definitions, $(\mathcal L, \bullet, \Delta,
\mathbf 1, \epsilon, \delta)$ is a differential graded commutative Hopf algebra
with unit.
\eth

\begin{proof}
The Leibnitz rule and the graded commutativity of the shuffle product
hold just as for $\mathcal D$, and one can easily verify the 
coassociativity of the coproduct and the compatibility between product 
and coproduct. 
%
%
The coboundary operator 
$\delta$ induces a coboundary operator on $\mathcal L \otimes 
\mathcal L$ by setting $\delta(\Gamma'\otimes\Gamma'') = 
\delta(\Gamma')\otimes \Gamma'' +  (-1)^{|\Gamma'|} \Gamma' \otimes 
\delta(\Gamma'')$ on the generators and extending this definition 
to $\mathcal L \otimes \mathcal L$ by linearity. Then one has the 
following equation
\[
\Delta\, \delta = \delta\,\Delta.
\]
In fact, it easy to check that, for any graph $\Gamma$, the explicit 
expression for $\Delta(\delta \Gamma)$ contains the same graphs (with the 
same signs) of the explicit expression for $\delta\Delta(\Gamma)$.
%
%
Finally, since this bialgebra is $\N$-graded and the only element in degree 
zero is the unit, the antipode is automatically defined.
%
%
%
\end{proof}

\begin{Rem}
If we pass to the dual, then our algebraic structures 
correspond to the ones considered in \cite{T} (see also \cite{BN95}
for the degree zero in the odd case). 
\end{Rem}

\section{Trivalent cocycles of even type}
\label{sec:inf}

A consequence of the algebra structure on $H^{k,0}(\mathcal D_e)$ is 
the following 

\bprop
\label{nonzerop}
For every $k\in\N$ we have $H^{2k,0}(\mathcal D_e)\neq 0$
\eprop

It is well known that the odd-case analogue of this Proposition holds, 
thanks to the existence of Vassiliev knot invariants at any 
order \cite{BN95}.

\begin{proof}
Let us denote by the $\Psi = \frac14\Phi - \frac13 \Gamma_1$ 
the graph cocycle of figure~\ref{ec2}. 
\begin{figure}[h] 
\resizebox{6cm}{!}{\includegraphics{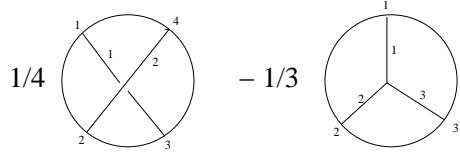}} 
\caption{Even cocycle of order 2} \label{ec2}
\end{figure}
We want to show that, for every $l\geq 2$, the $l$th power of $\Psi$ is 
nontrivial. First, we notice that $\delta_e(\Psi)=0$ and hence, thanks to 
Proposition~\ref{prop:db}, we also have $\delta_e(\Psi^{\bullet l})=0$. 
Moreover $\Psi^{\bullet l}$ cannot be $\delta_e$-exact since it has degree
zero and there are no graphs of negative degree. Therefore, it is 
enough to prove that, for every $l\ge 2$, the cocycle $\Psi^{\bullet l}$ is 
a linear combination of graphs with at least one coefficient different from 
zero. 

Let us denoted by $\Gamma_l$ the graph represented in figure~\ref{xyn}. This 
graph has $l$ triples of edges are attached on the external 
circle in such a way that they do not overlap.  
 
\begin{figure}[h] 
\resizebox{4cm}{!}{\includegraphics{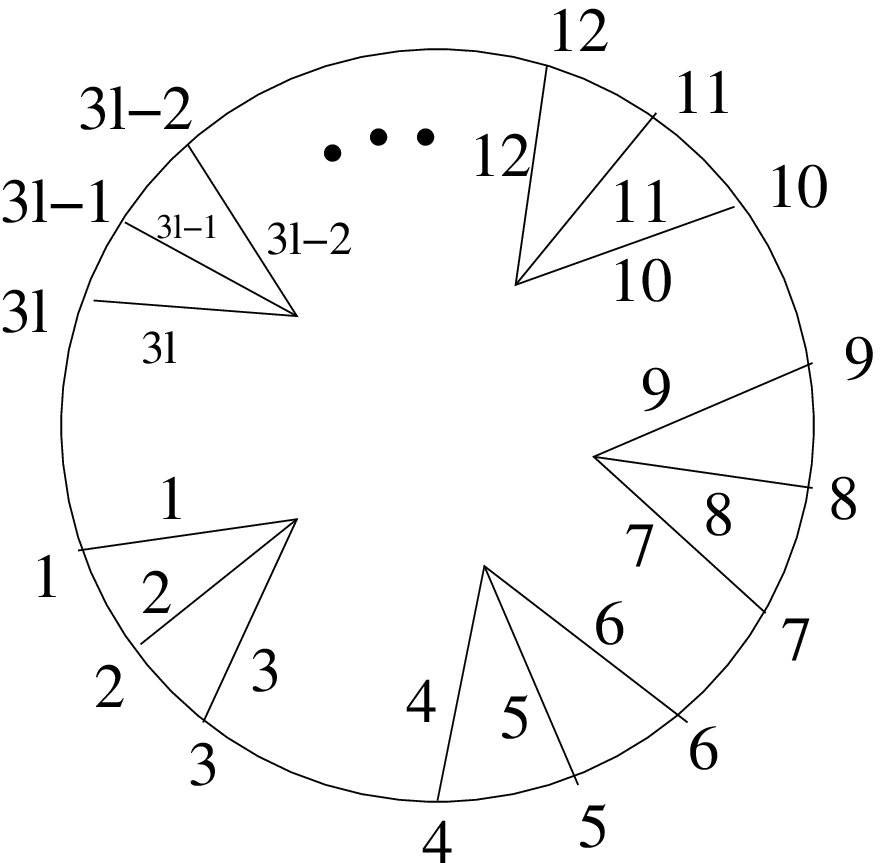}} \caption{} 
\label{xyn} 
\end{figure} 

We will prove by induction that the coefficient of $\Gamma_l$ in 
$\Psi^{\bullet l}$ is different from zero for every $l\geq 2$. The
proof of this fact for $l=2$ goes as follows: in the shuffle product 
\[
\Psi\bullet\Psi = \frac1{16}\ \Phi\bullet\Phi - \frac16\  \Phi\bullet\Gamma_1
+ \frac19\  \Gamma_1\bullet\Gamma_1
\]
there are exactly three contributions to 
$\Gamma_2$, all coming from $\Gamma_1\bullet\Gamma_1$. Since each 
of these three contributions has coefficient $\pm 1$, we have that the 
coefficient of $\Gamma_2$ in $\Psi^{\bullet 2}$ is different from zero.
 
Next, we suppose that the coefficient of $\Gamma_{l-1}$ in 
$\Psi^{\bullet (l-1)}$ is different from zero, and we observe that  
all the contributions to the graph $\Gamma_l$ in 
$\Psi^{\bullet l}$ arise 
from the shuffle product of $\Gamma_{l-1}$ with $\Gamma_1$. An easy 
computation shows that these contributions are all equal up to an even 
permutation of labels, and hence they cannot cancel out.
\end{proof} 

It is known from \cite{CCRL}, Theorem 1.1, that the maps 
\[
I\colon H^{k,0}(\mathcal D)\to H^{(n-3)k}(\imbr n)
\]
are injective for every $n>3$. This means that each of the graph 
cocycles of the above Proposition produces a nontrivial cohomology 
class of $\imbr n$ for every even $n\ge 4$. Hence we have: 
\bcor
For any $n> 3$ and for any positive integer $k_0$,
there are nontrivial cohomology classes on $\imbr{n}$
of degree greater than $k_0$.
\ecor

\section{Integration map}
\label{sec:imap}

We now recall how the maps of eqs.~\eqref{chainmap1} and \eqref{chainmap2} 
are constructed (see \cite{BT, CCRL} for further details). 
We consider the fiber bundle $p\colon C_{q,t}(\R^n) \to \imbr n$ 
whose fiber over a given imbedding $\gamma$ is the compactified 
configuration space of $q+t$ points in $\R^n$, the first $q$ of which are
constrained on $\gamma$. 

Let us fix a symmetric volume form $\omega^{n-1}$ on 
$S^{n-1}$, namely a normalized top form satisfying the additional condition 
$\alpha^*\omega^{n-1}=(-1)^n \omega^{n-1}$, where $\alpha$ is the antipodal 
map. A tautological form $\theta_{ij}$ is by definition the pull-back to 
$C_{q,t}(\R^n)$ of $\omega^{n-1}$ via the smooth map 
$\phi_{ij}\colon C_{q,t}(\R^n) 
\to S^{n-1}$ which, on the interior of $C_{q,t}(\R^n)$, is defined as 
\[
\phi_{ij} (x_1,\ldots,x_{q+t})=\frac{(x_i-x_j)}{|x_i-x_j|}\ . 
\]
For $i=j$ we use instead the map:
\[
C_{q,t}(\R^n)\stackrel{\pi}\lora C_{q,0}(\R^n)=C_{q}\times \imbr n 
\stackrel{pr_i\times id}\lora S^1\times \imbr n \stackrel{D}{\lora}
S^{n-1}
\]  
where $C_q$ is a component of the compactified configuration space of
$q$ points on $S^1$, $\pi$ forgets the $t$ points not lying on the imbedding,
$pr_i$ is the projection on the $i$th point and $D$ is the normalized
derivative $D(t,\psi)=\dot\psi(t)/|\dot\psi(t)|$.

For any given graph $\Gamma\in\mathcal D^{k,m}$ with $q$ external 
vertices and $t$ internal vertices, we construct a differential form
$\omega(\Gamma)$ on $C_{q,t}(\R^n)$ by associating the tautological form 
$\theta_{ij}$ to the edge connecting the vertices $i$ and $j$, 
and taking the wedge product of these forms over all the edges of $\Gamma$.
Then $I(\Gamma)$ is set to be the integral of $\omega(\Gamma)$ 
along the fibers of $p\colon C_{q,t}(\R^n) \to \imbr n$. 
The map $I$ extended by linearity to 
$\mathcal D^{k,m}$ takes value in $\Omega^{(n-3)k+m}
(\imbr n)$. Similarly, one can define $I\colon \mathcal L^{k,m} \to 
\Omega^{(n-3)k+m}(\imbrr n)$.\\ 

Let us now turn to the algebraic structures in the case of 
imbeddings of $S^1$ into $\R^n$. We know from \cite{CCRL}, Theorem 4.4, 
that the integration map $I\colon \mathcal{D}\to
\Omega(\imbr n)$ is a chain map for $n>3$. 
An easy check shows that the shuffle product
correspond exactly to the wedge product of configuration
space integrals. Therefore we have:

\bprop
The integration map $I\colon \mathcal{D}\to\Omega(\imbr n)$ is
a homomorphism of differential algebras with unit for $n>3$.
\eprop

Next we consider $\imbrr n$. 
\blem
\label{chl}
The integration map $I\colon \mathcal{L}\to\Omega(\imbrr n)$ is
a chain a map for $n>3$.
\elem

\begin{proof}
The only difference with respect to the case of imbeddings of the circle
is that now one has to consider also faces describing points on the 
imbedding escaping to infinity, possibly along with external points. 
A main feature of the compactified configuration spaces is that they 
split near the codimension-one faces in a product of the configurations 
collapsing to a certain point (or escaping at the point ``infinity'') and 
the configurations which remain far from this collapsing point. 
Since the imbeddings are fixed outside a compact set, the integration
along the faces at infinity yields zero unless the form degree of the 
integrand is zero on the first component (i.e., the configurations escaping 
at infinity). This implies that an entire 
connected component of a graph has to escape to infinity:
in fact, whenever exactly one argument of a tautological form goes 
to infinity, the form degree is entirely carried by the point escaping at 
infinity as a consequence of
\[
\frac{x-y}{|x-y|} \sim \frac x{|x|}\qquad \text{for $x\to\infty$}.
\]
If however a connected subgraph $\Gamma$ yields a zero form after 
integration on the face at infinity, the relation $(n-1)e=v_e+nv_i-1$ 
should hold (where again $e$ is the number of edges, $v_e$ the number 
of external vertices and $v_i$ the number of internal vertices of $\Gamma$). 
But in a graph whose vertices are at least trivalent $(n-1)e-v_e-nv_i$ is 
nonnegative.
\end{proof}

The proof of the next Proposition is exactly as for $\imbr n$.

\bprop
The integration map $I\colon \mathcal{L}\to \Omega(\imbrr n)$ 
is a homomorphism of differential algebras with unit for $n>3$.
\eprop

\section{Hopf algebras}
\label{sec:hopf}

We now want to define a coproduct on $\Omega(\imbrr n)$. In the
following we will call {\em long knots} the elements of $\imbrr n$.
We first fix some convention in the
definition of $\imbrr n$; namely, we choose a basis $\{\mathbf{e}_1,\dots,
\mathbf{e}_n\}$ of $\R^n$ such that all the elements of
$\imbrr n$ coincide, outside a compact subset,
with the reference imbedding $\varsigma(t)=\mathbf{e}_1\,t$.
Then we observe that two long knots $\gamma_1$ and $\gamma_2$ can be composed
to a new long knot $m(\gamma_1,\gamma_2)$ by
\[
m(\gamma_1,\gamma_2)(t)=
\begin{cases}
\Phi_-(\gamma_1(\phi_-(t)))  & t\le0,\\
\Phi_+(\gamma_2(\phi_+(t))) & t\ge0.
\end{cases}
\]
where $\Phi_\pm(x_1,x_2,\dots,x_n) = (\phi_\pm^{-1}(x_1),x_2,\dots,x_n)$
and $\phi_\pm\colon\reali^\pm\to\reali$ are any pair of diffeomorphisms, 
e.g., $\phi_\pm(t) = \tan\left(2\arctan (t) \mp\frac\pi2\right)$. Roughly
speaking we are attaching  $\gamma_1$ and $\gamma_2$ one after the other. 

The product $m$ has a unit $e$ given by the linear imbedding $\varsigma$ 
used to define $\imbrr n$.
The pullback of the product $m\colon\imbrr n\times \imbrr n\to\imbrr n$
defines a non-coassociative coproduct
\[
m^*\colon\Omega(\imbrr n)\to\Omega(\imbrr n)\hat\otimes
\Omega(\imbrr n).
\]
where $\hat\otimes$ denotes the topological tensor product of Fr\'echet 
spaces. The evaluation at $e$ defines a counit $\eta$:
\[
\eta(\omega)=
\begin{cases}
\omega(e) & \text{if $\deg\omega=0$}\\
0 & \text{if $\deg\omega>0$}.
\end{cases}
\]

We observe however that the non-coassociative coproduct $m^*$ gives rise to an 
associative operation in cohomology (in fact the product $m$ is associative up 
to a homotopy given by composing the long knot with a suitable diffeomorphism 
of $\R$). More precisely, we have:

\bprop
$(H(\imbrr n),\wedge,1,\Delta,\eta)$
is a graded commutative and cocommutative Hopf algebra for $n>3$.
\eprop

\begin{proof}
As shown in \cite{B}, there exists an action of the little 2-cubes operad 
$\mathcal {LC}_2$ on $\imbrr n$. This means that there are operations on the 
space $\imbrr n$ corresponding to each element of $\mathcal {LC}_2$. 
In particular, it turns out that one of these operations is the 
multiplication $m$ described above. Passing on the cochain level we obtain 
an action of operad $\Omega(\mathcal {LC}_2)$ on $\Omega(\imbrr n)$, which 
give rise in cohomology to an action of the operad $H(\mathcal {LC}_2)$ on 
$H(\imbrr n)$. Since the operad $H(\mathcal {LC}_2)$ is the linear dual 
of the Gerstenhaber operad, we have that the coproduct is 
coassociative and cocommutative. The compatibility between the wedge product 
and the coproduct is obvious. Finally,
the existence of the antipode follows from the fact that $H(\imbrr n)$ is 
$\N$-graded with only one element in degree zero.
\end{proof}

\begin{Rem}
A more explicitly proof of the cocommutativity of the coproduct is 
based on the observation that the composition $m$ of two long knots is 
commutative up to homotopy. In fact, one can shrink one of the two long 
knots in a very small region and slide it along the other long knot (see 
\cite{T, Tth} for details). 
\end{Rem}

Our central result is then the following
\bth
\label{main}
The map $I\colon H(\mathcal{L})\to H(\imbrr n)$ is
a Hopf algebra homomorphism for $n>3$.
\eth

\begin{proof}
First, we notice that for $n>3$, the degree of the differential form 
$I(\Gamma)$ is zero if and only if $\Gamma = \mathbf 1$. 
Therefore we have $\eta(I(\Gamma))=\epsilon(\Gamma)$.

We now have to show that the coproducts are compatible; viz.,
\begin{equation}
\label{copro}
\Delta\circ I = (I\otimes I)\circ\Delta.
\end{equation}
To prove this, consider any two cycles $Z_1$ and $Z_2$ (of degree
$k$ and $l$) of imbeddings. Let $Z=m_*(Z_1,Z_2)$ be the
$(k+l)$-cycle obtained by attaching the two cycles:
\[
Z(u_1,\dots,u_k,v_1,\dots,v_l)(t)
=
\begin{cases}
\Phi_-(Z_1(u_1,\dots,u_k)(\phi_-(t)))  & t\le0,\\
\Phi_+(Z_2(v_1,\dots,v_l)(\phi_+(t))) & t\ge0.
\end{cases}
\]

Identity~\eqref{copro} is equivalent to 
\begin{equation}\label{cigamma}
\int_Z I(\Gamma) = \sum (-1)^{\lambda(\Gamma',\Gamma'')}
\int_{Z_1} I(\Gamma')\,\int_{Z_2} I(\Gamma''),
\end{equation}
for any $Z_1$ and $Z_2$ and for any $\Gamma\in H(\mathcal{L})$,
where we write $\Delta\Gamma=\sum(-1)^{\lambda(\Gamma',\Gamma'')}
\Gamma'\otimes\Gamma''$.

To prove \eqref{cigamma}, let us introduce the $(k+l+1)$-chain $\mathfrak Z$ 
by
\begin{multline}
\label{chain-zeta}
\mathfrak Z(R,u_1,\dots,u_k,v_1,\dots,v_l)(t)
=\\
=\begin{cases}
\Phi_-(Z_1(u_1,\dots,u_k)(\phi_-(t+R)))-\varsigma(R)  & t\le-R,\\
\varsigma(t) & -R<t<R,\\
\Phi_+(Z_2(v_1,\dots,v_l)(\phi_+(t-R)))+\varsigma(R) & t\ge R,
\end{cases}
\end{multline}
with $R\in[0,+\infty)$. In practice we are moving the support of the cycles 
$Z_1$ and $Z_2$ far apart, and the parameter $R$ measure the distance between 
the two cycles. Since $\delta\Gamma=0$ and since $I$ is a chain map 
(Lemma~\ref{chl}), by Stokes' Theorem we get
\[
0=\int_{\partial \mathfrak Z} I(\Gamma) =
-\int_Z I(\Gamma) + \lim_{R\to+\infty}J_R
\]
with
\[
J_R(\Gamma) = \int_{\mathfrak Z(R,\cdot)} I(\Gamma).
\]
Let us write our graph cocycle as $\Gamma=\sum_i c_i\Gamma_i$. Each
$I(\Gamma_i)$ can be split into two parts as follows.
We fix $R\in (0,+\infty)$ and suppose $\Gamma_i$ has $q$ external vertices 
and $t$ internal vertices. Then we define $C^{1,R}_{q,t}(\R^n)$  to be the
subbundle of $C_{q,t}(\R^n)$ where the points corresponding to the 
external vertices of every primitive subgraph of $\Gamma_i$ lie either 
all to the left of $R$ or all to the right of $-R$. We also let 
$C^{2,R}_{q,t}(\R^n)$ be the fiberwise complement of 
$C^{1,R}_{q,t}(\R^n)$ in 
$C_{q,t}(\R^n)$. Now we define $I_{\alpha,R}(\Gamma_i)$, $\alpha=1,2$, 
to be the integral of the differential form $\omega(\Gamma_i)$ performed 
along the fibers of $C^{\alpha,R}_{q,t}(\R^n)$. Finally we set 
$I_{\alpha,R}(\Gamma) = \sum_i c_i I_{\alpha,R}(\Gamma_i)$ and we have
$J_R(\Gamma) = \int_{\mathfrak Z(R,\cdot)} I_{1,R}(\Gamma)  + 
\int_{\mathfrak Z(R,\cdot)} I_{2,R}(\Gamma)$. 

One immediately sees that $\lim_{R\to\infty} \int_{\mathfrak Z(R,\cdot)}
I_{1,R}(\Gamma)$ is equal to the right-hand side of eq.~\eqref{cigamma}, and 
hence what we have to prove is that 
\begin{equation}
\label{schif}
\lim_{R\to\infty} 
\int_{\mathfrak Z(R,\cdot)} I_{2,R}(\Gamma_i)=0.
\end{equation} 

We now need a generalization of Lemma 10 of \cite{AF}. Let $\Gamma$ be a 
connected graph with $v_e$ external vertices and $v_i$ internal vertices, 
and let $\omega(\Gamma)$ be the product of tautological forms associated to 
the graph $\Gamma$. We denote by $g_\Gamma$ 
the integral of $\omega(\Gamma)$ over the internal vertices of $\Gamma$, 
namely the push-forward of $\omega(\Gamma)$ along the map 
$p\colon C_{v_e,v_i}(\R^n)\to C_{v_e,0}(\R^n)$ that forgets the internal 
points. Let us write $g_\Gamma$ in coordinates as $g_\Gamma=
g_{\Gamma_I}(x_1,\ldots, x_{v_e}, \gamma)\, \dd x^I$, where $I$ is a 
multi-index. Suppose moreover that $\mathbf{x}(T) = 
(x_1(T),\ldots, x_{v_e}(T))$ a sequence in the configuration space with the 
property that there is a pair of points whose distance diverges as $T$ goes 
to infinity.

\blem
\label{lemma1}
With the above notations, we have
\[
\lim_{T\to\infty} \left. g_{\Gamma}\right|_{\mathbf{x}(T)} =0.
\]
\elem

\begin{proof} {\sf Case 1.} 
We consider first the case when $\Gamma$ has no edges whose end-points 
are both external. Then $g_\Gamma$ turns out to be bounded whenever its 
arguments $x_1,\ldots, x_{v_e}$ run in a bounded subdomain of $\R^n$ 
(namely, $g_\Gamma$ do not diverge if two or more points collapse). We 
can always suppose that the pair of points whose distance diverges 
fastest are $x_1(T)$ and $x_{v_e}(T)$. Using the translation invariance of 
the integral, we also suppose $x_{v_e}(T)=0$. Therefore we have that 
$|x_1(T)|\to\infty$ as $T\to\infty$.

We claim that if we rescale all the variables by $1/|x_1|$ we get
\[
g_{\Gamma_I}(x_1, \ldots, x_{v_e-1},0) = 
\left(\frac1{|x_1|}\right)^{\alpha_\Gamma} 
g_{\Gamma_I}(x_1/|x_1|, \ldots, 
x_{v_e-1}/|x_1|,0).
\]
where $\alpha_\Gamma = (n-1)e -n v_i$. In fact, when performing the 
integral along the fibers of $p\colon C_{v_e,v_i}(\R^n)\to C_{v_e,0}(\R^n)$, 
it is convenient to rescale the integration variables by $1/|x_1|$. 
This yields the contribution $-n v_i$. At this point the tautological forms 
$\theta_{ij}$, whose degree is $n-1$ and whose number is the same as the 
number of edges $e$, are obtained by rescaling the arguments of the 
functions $\phi_{ij}$, and this yields the contribution $(n-1)e$. 
Moreover, using the fact that $2e-3v_i-v_e\geq 0$ and $v_e>0$, we have 
\[
\alpha_\Gamma = (n-1)e -n v_i \geq (n-1)\frac32v_i + \frac{n-1}2v_e - nv_i =
\frac{n-3}2v_i + \frac{n-1}2v_e >0. 
\]
Now, $x_k/|x_1|$ is bounded for every $k=2,\ldots,v_e-1$, and hence 
the quantity $g_{\Gamma_I}(x_1/|x_1|, \ldots, x_{v_e-1}/|x_1|,0)$
remains bounded when $|x_1|$ goes to infinity. This proves the first case.

{\sf Case 2.}
When $\Gamma$ has edges whose end-points are both external, we 
denote by $\widetilde\Gamma$ the same graph with these edges removed.
We have then $g_\Gamma = g_{\widetilde\Gamma}h_\Gamma$ where 
$h_\Gamma$ is the 
product of the functions associated to the removed edges.
If $\widetilde\Gamma$ is the empty graph, then there is at least 
one edge whose end-points go far apart. Since the imbeddings are fixed 
outside a compact set, asymptotically $h_\Gamma$ is the pull-back of a volume 
form via a constant function, i.e., $h_\Gamma$ vanishes in the limit 
$T\to\infty$. If $\widetilde\Gamma$ is non empty but connected we are done 
by the argument in Case 1. The last possibility 
is when $\widetilde\Gamma$ is non empty and
disconnected. If inside a connected component there is a pair of points 
whose distance diverges, we are also done by the argument in Case 1, otherwise
there is at least one edge whose end-points go far apart, and in this case 
$h_\Gamma$ goes to zero.
\end{proof}

We now show that eq.~\eqref{schif} holds. 
Suppose $\Gamma_i$ has $q$ external vertices and $t$ internal vertices, and 
consider $C^{2,R}_{q,t}(\R^n)$.
For every element of $C^{2,R}_{q,t}(\R^n)$ consider a primitive 
subgraph $\Xi$ of $\Gamma_i$ such that the preimage $l$ of its 
left-most point on the imbedding 
is less than $-R$, while the preimage $r$ of its right-most point is 
greater than $R$. 
Let $\Xi_0, \ldots, \Xi_k$ be the connected component of 
$\Xi$ minus the oriented line, ordered following the order of their 
left-most external vertices. With the symbols $l_i$ and $r_i$ we mean 
the preimages of the left-most and right-most among the points 
corresponding to the external vertices of $\Xi_i$. In particular, 
$l_0 = l$. If $r_0>0$ then the result follows from Lemma~\ref{lemma1}
applied to $g_{\Xi_0}$. If $r_0<0$ then $l_1<0$. If $l_1<-R$ 
then we repeat our considerations for $\Xi/\Xi_0$, otherwise we have 
two possibilities: $r_1<R$, and the integral of $\omega(\Gamma_i)$ on 
this portion of $C^{2,R}_{q,t}(\R^n)$ is zero for 
dimensional reasons, or $r_1>R$ and 
we can apply Lemma~\ref{lemma1} to $g_{\Xi_1}$.

This concludes the proof of the Theorem.
\end{proof}

\section{The ``anomaly'' term}
\label{sec:anomaly}

Let us consider now the case of ordinary framed long knots, i.e., 
imbeddings of $\R$ into $\R^3$ with a choice of a trivialization of  
their normal bundle. Following \cite{BT}, we know that $I$ is not a 
chain map because of an anomaly term due to contributions from the 
most hidden faces of the boundaries of the compactified configuration 
spaces. More precisely, if we restrict to trivalent graphs, the following 
equation holds 
\[
\dd \tilde I = I\delta
\]
where 
\[
\label{i-tilde}
\tilde I(\Gamma)=I(\Gamma)+c_\Gamma I(\Theta).
\]
Here $c_\Gamma$ are certain (unknown) coefficients while $\Theta$ is 
the graph with one chord only, so $\int_K I(\Theta)$ is the self-linking 
number $\mathrm{slk}(K)$ of the framed long knot $K$. 
As a consequence, the Bott--Taubes map $\tilde I$ associates to each 
cocycle of trivalent graphs an invariant of framed long knots.

Let $\Gamma$ be a cocycle of trivalent graphs, $Z_1$ and $Z_2$ two
framed long knots, $Z=m(Z_1,Z_2)$ their product and $\mathfrak Z$ 
the 1-chain of eq.~\eqref{chain-zeta}. Then $\dd I(\Gamma) +
c_\Gamma \dd I(\Theta) = \dd \tilde I (\Gamma)=   I
(\delta(\Gamma)) =0$ and hence, applying Stokes' Theorem, we have 
\[
\lim_{R\to\infty}\int_{\mathfrak Z(R)}I(\Gamma)  - \int_Z I(\Gamma) = 
\int_{\mathfrak Z} \dd I(\Gamma) = 
 - c_\Gamma (\mathrm{slk}(\mathfrak Z(R)) - \mathrm{slk}(Z)) = 0,
\]
since the self-linking number of $\mathfrak Z(R)$ does not 
depend on $R$ and it is therefore equal to the self-linking number of $Z$. 
As a consequence, using the same arguments of the proof of  
Theorem~\ref{main}, we deduce that also in three dimensions 
the map $I$, when restricted to trivalent  
graphs, is compatible with the coproducts. On the other hand the 
Bott--Taubes map $\tilde I$ is a map of coalgebras if and only if 
\begin{equation}
\label{Itilde} 
\int_{Z} \tilde I(\Gamma) = \sum (-1)^{\lambda(\Gamma',\Gamma'')}
\int_{Z_1} \tilde I(\Gamma')
\int_{Z_2} \tilde I(\Gamma'').  
\end{equation}
The left-hand side of \eqref{Itilde} is given by 
\begin{multline*} 
\int_Z I(\Gamma) + c_\Gamma\int_Z I(\Theta) =  \sum \int_{Z_1} I(\Gamma') 
\int_{Z_2} I(\Gamma'') + c_\Gamma \left(\mathrm{slk}(Z_1) +\mathrm{slk}(Z_2)  
\right), 
\end{multline*} 
while the right-hand side is 
\begin{multline*}
 \sum\int_{Z_1} (I(\Gamma') + c_{\Gamma'}I(\Theta))
 \int_{Z_2}(I(\Gamma'') + c_{\Gamma''}I(\Theta)) =\\
= \sum\int_{Z_1}I(\Gamma')\int_{Z_2} I(\Gamma'') +
 c_{\Gamma''} \mathrm{slk}(Z_2) \int_{Z_1}I(\Gamma') + c_{\Gamma'}
 \mathrm{slk}(Z_1)\int_{Z_2} I(\Gamma'') +\\ +  c_{\Gamma'}c_{\Gamma''} 
\mathrm{slk}(Z_1) \mathrm{slk}(Z_2)
\end{multline*}
In particular we can write eq.~\eqref{Itilde} in the case when $\Gamma$ is 
obtained by ``attaching'' the oriented lines of two graphs $\Gamma_1$ 
and $\Gamma_2$. By using the fact of \cite{AF} that $c_\Gamma=0$ if $\Gamma$ 
is not primitive, then one easily sees that \eqref{Itilde} holds if and only 
if 
\[
c_{\Gamma_2}\mathrm{slk}(Z_2)\int_{Z_1} I(\Gamma_1) + 
c_{\Gamma_1}\mathrm{slk}(Z_1)\int_{Z_2} I(\Gamma_2) +
c_{\Gamma_1}c_{\Gamma_2}\mathrm{slk}(Z_1)\mathrm{slk}(Z_2)=0.
\]
Hence, if $c_{\Gamma_1}$ is different from zero, $\tilde I$ is not a
coalgebra map (e.g. choose $Z_1$, $Z_2$ and $\Gamma_2$ such that
$\mathrm{slk}(Z_2)=0$, $\mathrm{slk}(Z_1)\neq 0$ and
$\int_{Z_2}I(\Gamma_2) \neq 0$). In other words, the anomaly term
$c_\Gamma I(\Theta)$ is the obstruction for $\tilde I$ to be a
coalgebra map.

\end{document}